\documentclass{elsarticle}
\makeatletter
\def\ps@pprintTitle{%
 \let\@oddhead\@empty
 \let\@evenhead\@empty
 \def\@oddfoot{\centerline{\thepage}}%
 \let\@evenfoot\@oddfoot}
\makeatother
\usepackage{epsfig}
\usepackage{epic}
\usepackage{eepic}
\usepackage{enumerate}
\usepackage{graphics}
\usepackage{subfigure}
\usepackage{graphicx}

\usepackage[english]{babel}
\usepackage[utf8]{inputenc}
\usepackage{amsmath}
\usepackage{amssymb}
\usepackage{amsthm}
\usepackage{enumerate}
\usepackage{graphicx}
\usepackage[colorinlistoftodos]{todonotes}
\usepackage{breqn}

\newtheorem{thm}{Theorem}[section]
\newtheorem{lemma}[thm]{Lemma}

\newcommand{\sm}{\setminus}

\DeclareMathOperator{\po}{\overrightarrow{\chi}}

\begin{document}


\begin{frontmatter}

\title{Proper Orientations of Planar Bipartite Graphs}

\author[SFU]{Fiachra Knox\fnref{knox}}
\ead{fknox@sfu.ca}

\author[SFU]{Sebasti\'an Gonz\'alez Hermosillo de la Maza \corref{cor1}}

\ead{sga89@sfu.ca}

\author[SFU]{Bojan Mohar\fnref{mohar}}
\ead{mohar@sfu.ca}

\author[ParGO]{Cl\'audia Linhares Sales\fnref{sales}}
\ead{linhares@lia.ufc.br}

\address[SFU]{Simon Fraser University\\
						8888 University Drive\\
						Burnaby, BC, Canada\\}
						
\address[ParGO]{ParGO-DC, Federal University of Cear\'a\\
Cear\'a, Brazil}





	
\cortext[cor1]{Corresponding author}
 \fntext[knox]{Supported by a PIMS post-doctoral fellowship}
   \fntext[mohar]{Supported in part by an NSERC Discovery Grant (Canada), by the Canada Research Chair program, and by the Research Grant P1–0297 of ARRS (Slovenia)}
 
  \fntext[sales]{Partially supported by CAPES (Brazil) 99999.000458/2015-05}

\begin{abstract}
An orientation of a graph $G$ is {\it proper} if any two adjacent vertices have different indegrees. The {\it proper orientation number} $\po(G)$ of a graph $G$ is the minimum of the maximum indegree, taken over all proper orientations of $G$. 
In this paper, we show that a connected bipartite graph may be properly oriented even if we are only allowed to control the orientation of a specific set of edges, namely, the edges of a spanning tree and all the edges incident to one of its leaves. As a consequence of this result, we prove that 3-connected planar bipartite graphs have proper orientation number at most 6. Additionally, we give a short proof that $\po(G) \leq 4$, when $G$ is a tree and this proof leads to a polynomial-time algorithm to proper orient trees within this bound.
\end{abstract}


\end{frontmatter}


\section{Introduction}

Let $G=(V,E)$ be a simple graph, i.e., finite, with no loops or multiple edges. An {\it orientation}~$\sigma$ of $G$ is a digraph obtained from~$G$ by replacing each edge by exactly one of the two
possible arcs with the same endvertices. If $e=(u,v)$ is an arc of $\sigma$, we say that $u$ and $v$ are {\it tail} and {\it head} of $e$, respectively.
For each $v \in V(G)$, the \emph{indegree} of~$v$ in~$\sigma$, denoted  by~$d^-_{\sigma}(v)$, is the number of arcs with head~$v$. We denote by $\Delta^-(\sigma)$ the maximum indegree of $\sigma$.
We write $d^-(v)$ when the orientation is clear from the context.
An orientation $\sigma$ of $G$ is \emph{proper} if $d^-(u)\neq d^-(v)$, for all $uv\in E(G)$.
An orientation with $\Delta^-$ at most~$k$ is called a \emph{$k$-orientation}.
The \emph{proper orientation number} of a graph $G$, denoted by~$\po(G)$, is the minimum integer~$k$ such that~$G$ admits a proper $k$-orientation.

This graph parameter was introduced by Ahadi and Dehghan~\cite{AD13} in 2013.
They observed that this parameter is well-defined for any graph $G$ since one can always obtain a proper $\Delta(G)$-orientation. They also proved that deciding if a graph $G$ has proper orientation number equal to 2 is $NP$-complete even if $G$ is a planar graph.

Lately, Araujo et al.\ {\cite{AHLS16}} proved that the problem of determining the proper orientation number of a graph remains NP-hard for subclasses of planar graphs that are also bipartite and of bounded degree. In the same paper, they proved that the proper orientation number of any tree is at most 4. They raised the question about which classes of graphs have bounded proper orientation number, and particularly, if this would be the case for planar graphs. Recently, in an attempt to answer this question, two other classes of planar graphs were shown to have bounded proper orientation number, namely, cacti and claw-free planar graphs, with proper orientation number at most 7 and 6, respectively {\cite{AHLS16}}. These bounds are all tight. Despite of this progress, the main question remains open.

In this paper, we prove that a connected bipartite graph may be properly oriented even if we are only allowed to control the orientation of a specific set of edges: the edges of a spanning tree and all the edges incident to one of its leaves (see Lemma \ref{lem:orient bipartite}). As consequence of this result, we prove that 3-connected planar bipartite graphs have proper orientation number at most 6 (see Theorem \ref{thm:3-conn planar bipartite PO}). On the contrary of the other results on trees and cacti, our proof leads to a polynomial-time algorithm to find a proper orientation within this bound.

We close this paper with an alternative simpler proof of the bound for trees. Previous proofs did not give efficient way for finding a proper orientation with bounded indegrees. However, our proof results in a polynomial-time algorithm to properly orient this class of graphs within this bound.

\section{General bipartite graphs}

Determining the proper orientation number of regular and bipartite graphs has been a challenge ever since the parameter was defined. In the seminal paper on proper orientation number {\cite{AD13}}, it was proved that:

\begin{thm}[\cite{AD13}]
Let $G$ be a $(2k+1)$-regular graph. Then $\po(G) = k+1$ if and only $G$ is bipartite.
\end{thm}

Actually, if $G$ is any $k$-regular bipartite graph, then $\po(G) = \lceil \frac{k+1}{2} \rceil$. However, for general bipartite graphs, it turns out that determining the value of the proper orientation number is very difficult.

\begin{thm}[\cite{ACD+15}]
It is NP-complete to decide if $\po(G)= 3$ for a given planar bipartite graph $G$ with maximum degree 5.
\end{thm}

The next lemma shows that, in spite of the hardness of the problem even when restricted to bipartite graphs $G$, one can obtain a bound for $\po(G)$ that is in general better than the trivial bound $\Delta(G)$, and a corresponding proper orientation of $G$, both in polynomial time, even if we are only allowed to control the orientation of a small subset of edges of $G$.

\begin{lemma} \label{lem:orient bipartite}
Let $G$ be a connected bipartite graph on vertex classes $A$ and $B$ and let $x_0$ be a vertex in $A$.
Let $T$ be any spanning tree of $G$, in which $x_0$ is a leaf, and let $S$ be the set of edges of $G \sm T$ which contain $x_0$ .
Then for any orientation $\sigma'$ of $G \sm (S \cup T)$, there is an orientation $\sigma$ of $G$ which extends $\sigma'$, and in addition satisfies the following properties:
\begin{enumerate}[\rm (i)]
\item $d^-(x)$ is even for every $x \in B$;
\item $d^-(x)$ is odd for every $x \in A \sm \{x_0\}$;
\item $d^-(x_0) \leq 1$;
\item $\sigma$ is a proper orientation of $G$;
\item $\Delta^-(\sigma) \leq \Delta^-(\sigma') + \Delta(T) + 1$.
\end{enumerate}
\end{lemma}

\proof We begin by orienting every edge of $S$ away from $x_0$; together with our assumption that $x_0$ is a leaf of $T$, this ensures that (iii) holds.
For any orientation of $T$, let $A'$ be the set of vertices of $A \sm \{x_0\}$ which have even indegree and let $B'$ be the set of vertices of $B$ which have odd indegree.
To complete the orientation $\sigma$, we now orient $T$ in such a way that $|A'| + |B'|$ is minimized.

To see that $\sigma$ satisfies (i) and (ii), suppose for a contradiction that $A'$ is nonempty (the case when $B'$ is nonempty is similar).
Let $x \in A'$ and consider the unique path $P$ in $T$ from $x$ to $x_0$.
Form $T'$ from $T$ by reversing the orientation of every edge of $P$.
Then the parity of every vertex of $G$ except for $x$ and $x_0$ is unchanged, while the parities of $x$ and $x_0$ are changed.
Thus $|A'|$ is decreased by one and $|B'|$ is unchanged, contradicting the minimality of $T$.

To see that the orientation is proper, let us first observe that by (i) and (ii), the parities of $d^-(x)$ and $d^-(y)$ are different for any edge $xy$ of $G - x_0$.
Further, if $d^-(x_0) = 1$, then this is also true when $x = x_0$ or $y = x_0$.
Finally, if $d^-(x_0) = 0$, then $d^-(y) > 0$ for every vertex $y$ in the neighbourhood of $x_0$.
This proves (iv), and (v) is immediate from the definition of $\sigma$. 
\endproof

Now we apply this lemma to the case of $3$-connected bipartite planar graphs.
Firstly, we recall that every $3$-connected planar graph has a spanning tree of maximum degree at most $3$ {\cite{B66}}.

Next, observe from Euler's formula that any bipartite planar graph has maximum average degree less than $4$. Hence, it has an orientation with maximum indegree at most $2$ {\cite{H65}}.

We therefore obtain the following:

\begin{thm} \label{thm:3-conn planar bipartite PO}
Let $G$ be a $3$-connected planar bipartite graph.
Then $G$ admits a proper orientation with maximum indegree at most $6$.
\end{thm}

\proof Let $A$ and $B$ be the vertex classes of $G$ and let $T$ be a spanning tree of $G$ with maximum degree at most $3$.
Without loss of generality $A$ contains a leaf $x_0$ of $T$.
Let $S$ be the set of edges of $G \sm T$ which contain $x_0$.
Let $\sigma''$ be an orientation of $G$ with maximum indegree at most $2$, and let $\sigma'$ be the restriction of $\sigma''$ to $G \sm (S \cup T)$.
The theorem now follows immediately by applying Lemma~\ref{lem:orient bipartite}. 
\endproof

The proof of Lemma \ref{lem:orient bipartite} is constructive and yields a polynomial-time algorithm for finding the corresponding orientation $\sigma$. (A naive way gives time $O(n^2)$, but it is easy to do it in linear time.) The same can be said for Theorem \ref{thm:3-conn planar bipartite PO} since a spanning tree of maximum degree 3 (the proof of existence is based on discharging).

We observe also that a 2-orientation of a bipartite planar graph $G$ can be found in time $O(n^2)$.
Indeed, $G$ has a vertex $v$ of degree at most $3$; 
we can remove $v$, obtain a $2$-orientation of $G - v$ by induction and then add back $v$ with its edges oriented inwards.
If $d(v) < 3$ then we are done.
If $d(v) = 3$, let $S$ be the set of vertices of $G$ from which there is a directed path to $v$.
Then $G[S]$ has at most $2|S| - 4$ edges, and hence has a vertex $w$ of indegree at most $1$.
We can find such a $w$, along with a directed path from $w$ to $v$, in time $O(n)$ by exploring $S$ in the natural way.
We then reverse the orientation of each edge of the path to obtain a $2$-orientation of $G$.

\section{Trees}

The proper orientation number of trees was tackled in \cite{ACD+15} resulting in a tight constant upper bound for the parameter:

\begin{thm}[\cite{ACD+15}]
For any tree $T$, $\po(G) \leq 4$.
\end{thm}

The proof in \cite{ACD+15} is by contradiction and does not yield a polynomial-time algorithm to find a proper orientation for any tree within this bound. In what follows, we give an alternative and shorter proof of this theorem which results in a polynomial-time algorithm. We start by recalling a result of {\cite{AHLS16}} about proper orientation of paths.

\begin{lemma}[\cite{AHLS16}]\label{lem:ColoringPn}
Let $P=(v_1,\ldots,v_n)$ be a path on $n$ vertices, $n\neq 2$. Then, there exists a proper $2$-orientation of $P$ such that $v_1$ and $v_n$ have indegree 0.
\end{lemma}

We can use Lemma \ref{lem:ColoringPn} to observe that if $n \notin \{1,3\}$, then there exists a proper $2$-orientation of $P$ such that $v_1$ and $v_n$ have indegree 0 and 1, respectively. Additionally, if $n = 3$, then $P$ has two possible proper $2$-orientations where both extremities have degree 0 or 1.

\begin{thm}\label{thm:trees}
Let $T$ be a tree. Then $\po(G) \leq 4$.
\end{thm}

\proof
The proof is by induction on the order of $T$. Observe that the theorem is trivially true for any tree with at most 5 vertices. Now, suppose that $n\ge6$ and that the theorem is true for any tree with less than $n$ vertices, and consider $T$, a tree of order $n$. The claim is clear if $T$ is a path. We may therefore assume that $\Delta(T) \geq 3$. Hence $T$ contains an edge $e=xy$ such that $d(x) \geq 3$ and such that every vertex except $x$ in the connected component $U$ of $T-e$ that contains $x$ has degree at most~2.

By the induction hypothesis, $T-U$ has a proper 4-orientation $\sigma'$. We are going to extend $\sigma'$ to an orientation $\sigma$ of $T$,
such that $d^-_{\sigma}(v) = d^-_{\sigma'}(v)$, for every $v \in T-U$. In particular, we orient the edge $yx$ from $y$ towards $x$. For the rest, we consider the following three cases.

\smallskip

{\noindent {\bf Case 1:}} $d(x)  \geq 4$. Let $L(x) =\{3,4\} \setminus \{d^-_{\sigma'}(y)\}$.
Observe that $x$ is attached to at least 3 paths in $U$.  Then $\sigma'$ can be extended to $T$ in the following way. According to $L(x)$, orient two or three neighbors of $x$ in $U$ towards $x$. This will assure that $d^-_{\sigma}(x) \ne d^-_{\sigma}(y)$. Then we apply Lemma \ref{lem:ColoringPn} and the remark after the lemma to properly orient the remaining edges of $U$. Since the vertices of the paths attached to $x$ in $U$ have indegree at most 2, the obtained orientation $\sigma$ is proper with $d^-_{\sigma}(x) \in L(x)$.

\smallskip

{\noindent {\bf Case 2:}} $d(x) = 3$ and no neighbor of $x$ is a leaf. If $d^-_{\sigma'}(y) \neq 3$, then $\sigma'$ can be extended to $T$ by orienting all the neighbors of $x$ towards $x$, and by completing the orientation of the paths attached to $x$ in $U$ as claimed in Lemma \ref{lem:ColoringPn}.  Otherwise, $d^-_{\sigma'}(y) = 3$. Then, orient all the edges incident to $x$ (except for $xy$) out of $x$. Since no path is a leaf, the orientation of $T$ can be completed by making each neighbor of $x$ in the paths have degree 2, and this is always possible by the remark after Lemma \ref{lem:ColoringPn}.

\smallskip

{\noindent {\bf Case 3:}} $d(x) = 3$ and at least one of its neighbors in $U$ is a leaf. Let $w$ be that leaf vertex and let $P$ be the path starting at $x$ and using the remaining vertices in $U$. Then $\sigma'$ can be extended to $T$ as follows. We orient $P$ such that the edge incident with $x$ is oriented towards $x$ (Lemma \ref{lem:ColoringPn}). This will assure that the indegree of $x$ will be 2 or 3. Finally, we orient the third edge $xz$ incident with $x$ in such a way that $d^-_{\sigma}(x) \ne d^-_{\sigma}(y)$.
\endproof

\bibliographystyle{acm}
\bibliography{PO-bipartite-160715}

\end{document}